\documentclass[11pt,leqno]{amsart}
\usepackage{graphicx}
\usepackage{amssymb}
\usepackage{amsfonts}
\usepackage{amsmath}                        
\usepackage{latexsym}
\usepackage{amsthm}
\usepackage{amscd}

\newtheorem*{maintheorem}{Main Theorem}
\newtheorem{theorem}{Theorem}[section]
\newtheorem{lemma}[theorem]{Lemma}
\newtheorem{corollary}[theorem]{Corollary}
\newtheorem{proposition}[theorem]{Proposition}

\newtheorem{definition}[theorem]{Definition}

\newcommand\vol{\mathbf{volume}\,}
\newcommand\PD{\mathbf{PD}}
\newcommand\R{\mathbb R}
\newcommand\Rn{\mathbb R^n}
\newcommand\C{\mathbb C}

\newcommand\Val{\mathbf {Val}}
\newcommand\ValSO{\mathbf {Val}^{SO(n)}}
\newcommand\ValU{\mathbf {Val}^{U(n)}}

\newcommand\Hom{\mathbf{Hom}}

\newcommand\area{\mathbf{area}}

\newcommand \ot{\otimes}
\newcommand \PP{\mathbb P}

\begin{document}

\title{Structure of the unitary valuation algebra}
\author{Joseph H.G. Fu}
\thanks{Partially supported 
by NSF grant DMS-0204826.}
\address{Department of Mathematics \\
          University of Georgia\\ Athens, Georgia 30602  }
\email{fu@math.uga.edu}

\begin{abstract}
S. Alesker has shown that if $G$ is a compact subgroup of $O(n)$ acting transitively on the unit sphere $S^{n-1}$ then the vector space $\Val^G$ of continuous, translation-invariant, $G$-invariant convex valuations on $\R^n$ has the structure of a finite dimensional graded algebra over $\R$ satisfying Poincar\'e duality.  We show that the kinematic formulas for $G$ are determined by the product pairing.  Using this result we then show that the algebra $\Val^{U(n) }$ is isomorphic to $\R[s,t]/(f_{n+1}, f_{n+2})$, where $s,t$ have degrees 2 and 1 respectively, and the polynomial $f_i$ is the degree $i$ term of the power series $\log(1 + s +t)$.\end{abstract}

\maketitle

\section{Introduction}
In \cite{Hadwiger:1957}, Hadwiger showed that the vector space $\Val^{SO(n)}$ of continuous convex valuations on $\R^n$ invariant under the group $\overline{SO(n)}$ of orientation-preserving isometries has dimension $n+1$, with a basis consisting of the Minkowski ``Quermassintegrale", or intrinsic volumes in the terminology of \cite{Klain-Rota:2000}. An immediate consequence is the following form of the Principal Kinematic Formula of Blaschke: if $\Phi_0, \dots \Phi_n$ are the intrinsic volumes, indexed by degree, then there exist constants $c^k_{ij}$ such that 
$$
\int_{\overline{SO(n)}} \Phi_k(A \cap \bar g B) \, d\bar g = \sum_{i+j=n+k} c^k_{ij} \Phi_i(A) \Phi_j(B) 
$$
for all compact convex bodies $A,B \subset \Rn$ (cf. also \cite{Santalo:1978}). By applying the formula to appropriate lists of bodies $A,B$, one may then determine the constants $c_{ij}^k$ by explicit calculations of the integral. Though not essentially difficult, this procedure can be a bit troublesome, with many opportunities for computational errors; however, Nijenhuis \cite{Nijenhuis:1974} showed that if the basis $\{\Phi_i\}_{i=0}^n$ and the Haar measure $d \bar g$ are normalized appropriately then all of the $c_{ij}^k$ are equal to unity. He speculated that there might exist some underlying algebraic structure that would explain this fact.

More recently, in a series of fundamental papers \cite{Alesker:2001}, \cite{Alesker:2003}, \cite{Alesker:2004} S. Alesker has shown  that if $G$ is a compact subgroup of the orthogonal group $O(n)$ acting transitively on the unit sphere $S^{n-1}$, then the vector space $\Val^G$ of $G$-invariant translation-invariant continuous valuations carries the structure of a finite-dimensional commutative graded algebra over $\R$. Furthermore the resulting algebra satisfies Poincare duality: the top degree piece of $\Val^G$ is one-dimensional and occurs in degree $n$, and the pairing $\langle a,b\rangle :=$ the degree $n$ component of $ab$ is perfect. One of the results of the present article is to show that this algebra structure satisfies Nijenhuis's speculation, reducing the results of \cite{Nijenhuis:1974}, which originally appeared to be a kind of miracle obtained by laborious calculations, to an obvious triviality based on the simple structure of the algebra $\Val^{SO(n)}$.

Thus the case $G= SO(n)$ should be viewed as the ground case of a more subtle general theory that remains to be worked out in detail. The first serious case is the case $G= U(n)$. The main result of this paper is the explicit determination of the structure of $\Val^{U(n)}$ (Thm. \ref{thm:main} below):
\begin{maintheorem} The graded $\R$-algebra $\Val^{U(n)}$ is isomorphic to $\R[s,t]/(f_{n+1},f_{n+2})$, where the generators $s,t$ have degrees 2 and 1 respectively, and $f_j$ is the degree $j$ component of the power series $\log(1+s+t)$.
\end{maintheorem}
Because of Thm. \ref{thm:KG and PD} below, this result determines in principle the kinematic formulas for all of the $U(n)$-invariant valuations. Nevertheless the problem of writing them down explicitly remains open. In fact this is only one of several open problems arising from a comparison between the $U(n)$ and the $SO(n)$ theories, which we discuss these in the closing section.

{\bf Acknowledgements.} It is a pleasure to thank A. Abrams, S. Alesker, D. Benson, R. Howard, S. Mason, G. Matthews, D. Nakano, T. Shifrin and R. Varley for illuminating discussions at various stages of this work.

\section{General results} Throughout this section we let $V$ be a vector space over $\R$ of dimension $n < \infty$, endowed with a euclidean structure. Let $O(V)$ denote the corresponding orthogonal group, and fix a compact subgroup $G \subset O(V)$ that acts transitively on the unit sphere of $V$. Put $\mathcal K(V)$ to be the space of all compact convex subsets of $V$, endowed with the Hausdorff metric. If $r \in \R$, $x \in V$ and $K \in \mathcal K(V)$ then we put
\begin{align*}
x+ K &:= \{x+p: p \in K\},\\
rK&:= \{rp: p \in K\}.
\end{align*}

Denote by $\Val^G(V)$ or simply ($\Val^G$) the vector space of continuous functions $\phi:\mathcal K(V) \to \R$ enjoying the properties
\begin{itemize}
\item finite additivity: if $K,L, K \cup L \in \mathcal K(V)$ then
$$ \phi(K \cup L) = \phi(K) + \phi(L) -\phi(K\cap L); $$
\item translation-invariance: if $x\in V$ and $K \in \mathcal K(V)$ then $\phi(x+K) = \phi(K)$;
\item $G$-invariance: if $g \in G$ and $K \in \mathcal K(V)$ then
$$ g\phi(K):= \phi(g^{-1}K) = \phi(K).$$
\end{itemize}

We put
\begin{equation}
\omega_k:=  \frac{\pi^{\frac k 2}}{\Gamma(\frac{k+2} 2)}
\end{equation}
for the volume of the unit ball in $\R^k$.

\subsection{$\Val^G$ as an algebra}
We begin by listing some basic properties of $\Val^G$. Put $\overline G:= G \ltimes V$, equipped with a bi-invariant Haar measure. We leave the choice of normalizing constant undetermined for the moment.

 Recall that a valuation $\phi$ is said to have degree $i$ if $\phi(rK) = r^i \phi(K)$ for all $r \ge 0$. Put $\Val^G_i $ for the subspace of degree $i$ elements in $\Val^G$.
The first fact is a consequence of a theorem of P. McMullen \cite{McMullen:1977}.

\begin{theorem} 
\begin{equation}
\Val^G(V) = \bigoplus_{i=0}^n \Val^G_i(V).
\end{equation}
Furthermore $\Val_n^G(V)$ is one-dimensional, and is spanned by the volume.
\end{theorem}

The results listed next are simple consequences of results of Alesker (\cite{Alesker:2001},\cite{Alesker:2003},\cite{Alesker:2004}). In fact we give them only in a restricted form sufficient for the purposes of the present paper.
\begin{theorem} \label{thm:compendium}
\begin{itemize} \item $\dim_\R\Val^G < \infty.$
\item Given $K \in \mathcal K(V)$, define $\mu^G_K \in \Val^G$ by
\begin{align}
\mu_K^G(L)&:= \int_G \vol(L - g K) \, dg\\
&= \int_{\overline G} \chi(L \cap \bar g K)\, d\bar g.
\end{align}
Then there are $K_1,\dots,K_N\in \mathcal K(V)$ such that $\Val^G$ is spanned by $\mu_{K_1}^G,\dots,\mu_{K_N}^G$.
\item There is a natural continuous multiplication on the space of all continuous translation-invariant valuations that restricts to a multiplication on $\Val^G$, given as follows: if $\phi \in \Val^G$ and $K ,L\in \mathcal K(V)$, then
\begin{equation}
(\phi\cdot \mu_K^G)(L):= \int_{\overline G} \phi(L \cap \bar g K) \, d\bar g.
\end{equation}
Extending by linearity, the resulting product endows $\Val^G$ with the structure of a commutative graded algebra over $\R$, with unit element given by the Euler characteristic $\chi$.
\item Let $W \subset V$ be a linear subspace, and let $H \subset G$ be the stabilizer of $W$.  Suppose that $H$ acts transitively on the unit sphere of $W$. Then the natural restriction map $\Val^G(V) \to \Val^H(W)$ is a homomorphism of $\R$-algebras.
\item  The pairing $\Val^G\otimes \Val^G \to \Val_n^G \simeq \R$ given by
\begin{equation}\label{eq:pairing}
\mathbf{PD}:(a,b) \mapsto (ab)_n
\end{equation}
(degree $n$ piece of $ab$) is perfect.
\end{itemize}
\end{theorem}

Thus the pairing (\ref{eq:pairing}) may be thought of as a self-adjoint map ${\PD}:\Val^G \to (\Val^G)^*$. For the moment we leave unspecified the choice of linear isomorphism $\Val_n^G \sim\R$, and $\PD$ inherits this imprecision. It is trivial to see that $\PD$ is a homomorphism of $\Val^G$-modules, where the $\Val^G$-module structure on the dual space $(\Val^G)^*$ is 
\begin{equation}
(a \alpha)(b):= \alpha(ab),
\end{equation}
$a,b\in \Val^G, \alpha \in (\Val^G)^*$. Since the pairing is perfect, in fact $\PD$ is an isomorphism, and it is clearly graded in the sense that
 $\PD:\Val^G_i \to (\Val^G_{n-i})^*$.
In other words $\Val^G$ carries the structure of a graded Frobenius algebra (\cite{Curtis&Reiner}).

It is well known and trivial to prove:

\begin{lemma}\label{lm:uniqueness}
If $A$ is a finite-dimensional graded algebra over a field $k$, then any two graded $A$-module isomorphisms $A \to A^*$ differ by multiplication by a unit of $A$ of pure degree 0.
\end{lemma}

\subsection{Kinematic formulas}
Let $\mu_1,\dots,\mu_N$ be a basis for $\Val^G$. It is straightforward to deduce that given $\phi \in \Val^G$ there are constants $c_{ij}^\phi \in \R$ such that for all $K,L \in \mathcal K(V)$
\begin{equation}
\int_{\overline G} \phi(K \cap \bar g L) \, d\bar g  = \sum_{i,j= 1}^N c^\phi_{ij}\mu_i(K)\mu_i(L).
\end{equation}
This situation may be abbreviated by defining the map
$$k_G:\Val^G \to \Val^G \otimes \Val^G\simeq \Hom_\R((\Val^G)^*,\Val^G)$$
by 
$$ k_G(\phi) := \sum_{i,j} c^\phi_{ij} \mu_i \ot \mu_j.$$
Note that the precise definition of $k_G$ depends on the choice of normalization for the Haar measure $d\bar g$. For the time being we prefer to leave this unspecified.
It is straightforward to check that $k_G$ is a coassociative, cocommutative coproduct. Noticing the similarity with the definition of the product, coassociativity is equivalent to
\begin{lemma}\label{lm:sesqui}
If $\phi,\psi \in \Val^G$ then
$$ k_G(\psi\cdot \phi) =\sum_{i,j} c^\phi_{ij} (\psi\cdot\mu_i )\ot \mu_j = \sum_{i,j} c^\phi_{ij} \mu_i \ot (\psi\cdot\mu_j).$$
$\square$
\end{lemma}

\begin{proposition} For every $\varphi \in \Val^G$, 
$k_G(\varphi)$ is a homomorphism of $\Val^G$-modules when thought of as a map $(\Val^G)^*\to \Val^G$. Furthermore $k_G(1)$ is an isomorphism.
\end{proposition}
\begin{proof} That $k_G(\varphi)$ is a homomorphism of $\Val^G$-modules follows immediately from Lemma \ref{lm:sesqui}. 

To see that $k_G(1)$ is an isomorphism it is enough to prove surjectivity. By Thm. \ref{thm:compendium} above, any given valuation $\phi \in \Val^G$ may be written as 
\begin{align}
\phi&= \sum_{i=1}^N a_i \mu^G_{K_i}\\
&= \sum_{i=1}^N a_i\int_{\overline G} \chi(\ \cdot \ \cap \bar g K_i) \, d\bar g \\
& = \sum_{i=1}^N a_i\,k_G(1)(\ \cdot \ , K_i).
\end{align}
for some $a_1,\dots,a_N \in \R$. But this last expression is precisely the image under $k_G(1)$ of the element $\psi\mapsto \sum a_i \psi(K_i)$ of $(\Val^G)^*$.
\end{proof}

In view of Lemma \ref{lm:uniqueness} this gives

\begin{theorem}\label{thm:KG and PD} With appropriate choices of scaling factors,
$$k_G(1) = \PD^{-1} .$$
\end{theorem}

\subsection{The classical case: $G = SO(n)$}
The case of $G = SO(n)$ or $O(n)$ was substantially settled by Hadwiger in the 1950s. Nonetheless the perspective introduced by Alesker remains illuminating.

Hadwiger proved that there is exactly one $SO(n)$-invariant valuation on $\R^n$ in each degree between 0 and $n$, given by the coefficients of the polynomial giving the volume of a tubular neighborhood of variable radius $r$ (Steiner's formula). Thus the bodies $K_i$ above may be taken to be balls of $n+1$ distinct radii $r_0< r_1<\dots<r_n$.

An alternative approach is to take $K_i$ to be a disk of dimension $i$. Letting the radius of such a disk to tend to $\infty$ and normalizing appropriately, one arrives at the classical expression for the Hadwiger valuations in terms of intersections with affine subspaces:
if $\overline G(n,k)$ denotes the affine Grassmannian of $k$-planes in $\R^n$ then the valuation $\mu_i$ of degree $i$ may be expressed
$$ \mu_i(K) = \int_{\overline G(n,n-i)}\chi(K \cap \bar P) \, d\bar P.$$

Let $t:= \mu_1$. Then
\begin{equation}
t(K)= \lim_{r\to \infty}r^{1-n} \int_{\overline{SO(n)}} \chi(K \cap \bar g D_r^{n-1}) \, d\bar g,
\end{equation}
where $D_r^{n-1}$ is the disk of dimension $n-1$ and radius $r$. Therefore the definition of the product gives
\begin{align}
t^2(K) &= \lim_{r\to\infty} r^{2-2n} \int_{\overline{SO(n)}} \int_{\overline{SO(n)}} \chi(K \cap \bar g D_r^{n-1}\cap \bar h D_r^{n-1}) \, d\bar g \, d\bar h \\
& = \int_{\overline G(n,n-2)} \chi(K \cap \bar Q) \, d\bar Q\\
&= \mu_{2}(K).
\end{align}
Continuing in this way we arrive at the following result of Alesker:
\begin{theorem}\label{thm:valSOn}
$\Val^{SO(n)} \simeq \R[t]/(t^{n+1})$.
\end{theorem}

The Poincare duality pairing is obviously $\langle t^i, t^j\rangle = \delta_i^{n-j}$. As a map $\Val^{SO(n)}\to(\Val^{SO(n)})^* $ it takes $t^{n-i}$ to  $(t^i)^*$, where $1^*,t^*,\dots, (t^n)^*$ is the dual basis to the $t^i$. The principal kinematic formula $k_{SO(n)}(1) = \PD^{-1}$ thus takes $(t^i)^*$ to $t^{n-i}$, or in different terms
$$ k_{SO(n)}(1) = \sum_{i=0}^n t^i \ot t^{n-i}.$$
More generally, Lemma \ref{lm:sesqui} yields
$$ k_{SO(n)}(t^k) = \sum_{i+j = n + k} t^i \ot t^j.$$
Thus we recover a classical result of Nijenhuis \cite{Nijenhuis:1974}:
\begin{theorem}
There exists a graded basis $\mu_0,\dots,\mu_n$ for $\Val^{SO(n)}$ such that for an appropriate normalization of the Haar measure
$$ \int_{\overline{SO(n)}} \mu_k(K\cap \bar g L) \, d\bar g = 
\sum_{i +j = n + k} \mu_i(K) \mu_j(L),$$
i.e. the coefficients in the kinematic formulas for the $\mu_i$ in terms of the $\mu_i$ are all equal to unity.
\end{theorem}

\subsection{The orthogonal complement of $\Val^{SO(n)}$ in $\Val^G$} Returning to the case of a general group $G$ (transitive on the sphere of course), we take $V= \R^n$. Then $Val^{SO(n)} \subset \Val^G$. Put 
$$ \mathcal A^G_i:= \{\mu \in \Val^G_i: t^{n-i}\cdot\mu= 0\}$$
 and
$$\mathcal A^G := \bigoplus_{i=1}^{n-1} \mathcal A^G_i.$$
Thus $\mathcal A^G_i$ is a subspace of codimension 1 of $\Val^G_i$. Clearly
$\mathcal A^G_i = 0$ unless $1 \le i \le n-1$, and
\begin{equation}\label{eq:inclusion}
t\cdot \mathcal A^G_i \subset \mathcal A^G_{i+1}.
\end{equation}

(Actually Alesker has shown that $\Val^G_i$ is one-dimensional for $i = 1, n-1$, hence
$\mathcal A^G_1= \mathcal A^G_{n-1} = 0$ as well.)

\begin{proposition}\label{prop:annihilator} With appropriate scalings,
\begin{equation}\label{eq:annihilator}k_G(t^k) \equiv \sum_{i+j = n +k} t^i \ot t^j   \mod \mathcal A^G \ot \mathcal A^G, \ k = 0,\dots,n.\end{equation}
\end{proposition}
\begin{proof} The subspace $\mathcal A^G_i$ is the orthogonal complement of $t^{n-i}$ under the product pairing. Theorem \ref{thm:KG and PD} may be interpreted as saying that $k_G(1) $ is the associated pairing on the dual space, which implies (\ref{eq:annihilator}) for $k=0$. Now the general case follows from Lemma \ref{lm:sesqui} and the definition of $\mathcal A^G$.
\end{proof}

\begin{corollary}\label{cor:annihilator}
$$ \mathcal A^G = \{\mu \in \Val^G: \mu(D^n_r) = 0 \rm{ \ for \ all \ } r > 0\}.$$
\end{corollary}
\begin{proof}
  Let $\{\alpha_l\}_l$ be a basis for $\mathcal A^G$. Then
by Prop. \ref{prop:annihilator}, given $K \in \mathcal K(V)$
\begin{align}
\sum_{i +j =n}  t^i (K) t^j (D^n_r)+ \sum_{k,l} d_{kl} \alpha_k (K) \alpha_l(D^n_r ) &=\int_{\overline G} \chi(K \cap \bar g D^n_r)\, d\bar g\\
& = \int_{\overline {SO(n)}} \chi(K \cap \bar g D^n_r)\, d\bar g \\
&=  \sum_{i+j = n} t^i(K)t^j(D^n_r). 
\end{align}
Since the pairing $k_G(1)$ is symmetric  and nonsingular, so is its restriction to $(\Val^{SO(n)})^\perp = \mathcal A^G$, from which it follows that all $\alpha_l(D^n_r) = 0$. 
\end{proof}

{\bf Remark.} In fact the group $G$ plays no role here, as Alesker (\cite{Alesker:personal})has proved the following: 

{\sl Let $\Val_k(\R^n)$ denote the space of all degree $k$ translation-invariant degree $k$ valuations on $\R^n$.
Put 
$$\mathcal A_k:= \{\phi \in \Val_k(\R^n):t^{n-k}\cdot \phi = 0\}$$ 
and
$$\mathcal B_k:= \{\phi \in \Val_k(\R^n): \phi (D^n(r) ) = 0 {\rm \  for \  all  \ } r > 0\}.$$ 
 Then $\mathcal A_k = \mathcal B_k$.
}

\section{The unitary case} 
\subsection{Statement of the main theorem}

It is natural to consider next the case $V= \C^n$, $G= U(n)$. Let $\bar G_\C(n,k)$ denote the affine Grassmannian of {\it complex} $k$-planes in $\C^n$. We will assume that
$$\C^1\subset \C^2\subset \C^3\subset \dots$$
in the natural way, and will take the corresponding unitary groups and Grassmannians to be included in each other accordingly.
Denote the unit ball in $\C^k$ by $D_\C^k$. 

In order to state a precise result in this case we need to take more care in the choices of constants. We normalize the measures $d \bar g$ on $\overline{U(n)}$ and $d \bar h$ on $\overline{SO(2n)}$ so that, given any measurable set $E\subset \C^n$, 
\begin{equation}\label{eq:normalization}
d\bar g\left(\{\bar g\in \overline{U(n)}: \bar g (0) \in E\}\right) = d\bar h\left(\{\bar h\in\overline{SO(2n)}: \bar h (0) \in E\}\right)=  |E|,
\end{equation}
where $|E|$ is the Lebesgue measure of $E$. We take the measures on the real affine Grassmannians $\bar G(2n,k)$ and the complex affine Grassmannians $\bar G_\C (n, l)$ so that
\begin{align}
\label{eq:Gbar}
	d\bar P \left(\{\bar P \in \bar G(2n,k): \bar P \cap D_\C^n \ne \emptyset 		\}\right) &= \omega_{2n-k}, \\
\label{eq:GCbar}
	d\bar Q \left(\{\bar Q \in \bar G_\C(n,l): \bar Q \cap D_\C^n \ne 			\emptyset \}\right) &= \omega_{2n-2l} = \frac{\pi^{n-l}}{(n-l)!}. 
\end{align}

These measures are compatible with the measures on $\overline{SO(2n)}$ and $\overline{U(n)}$ in the 
following sense.  We define $V^{2n}_0,\dots, V^{2n}_{2n} \in \Val^{SO(2n)}$ by
\begin{equation}
V_i^{2n} (K): = \int_{\bar G(2n, 2n-i)} \chi (K \cap \bar P) \, d\bar P = d\bar P\left( \{\bar P: \bar P \cap K \ne \emptyset\}\right)
\end{equation}
and $W^n_0,\dots, W^n_n \in \Val^{U(n)}$ by
\begin{equation}
W^n_j(K):= \int_{\bar G_\C(n, n-j)} \chi (K \cap \bar Q) \, d\bar Q = d\bar Q\left( \{\bar Q: \bar Q \cap K \ne \emptyset\}\right),
\end{equation}
$K \in \mathcal K(\C^n)$. Then
\begin{align}
\label{eq:W disks}
	W^n_j(K) &=\lim_{R\to \infty}( \omega_{2n-2j} R^{2n-2j})^{-1} d\bar g (\{\bar g\in  \overline{U(n)} : 			\bar g D_\C^{2n-2j}(R) \cap K \ne \emptyset\}), \\
\label{eq:V disks}
	V_i^{2n}(K) &=\lim_{R\to \infty}( \omega_{2n-i} R^{2n-i})^{-1} d\bar h (\{\bar h \in  \overline{SO(2 n)}: 		\bar h D_\R^{2n-i}(R) \cap K \ne \emptyset\}).
\end{align}

Clearly
$$ \deg V^{2n}_i = i, \ \deg W^n_j = 2j.$$

Our main result is 
\begin{theorem}\label{thm:main}
\begin{itemize}
\item $\Val^{U(n)}$ is generated as an $\R$-algebra by  $V^{2n}_1$ and $W^n_1$.
\item Consider the graded polynomial algebra $\R[s,t]$, where $\deg t = 1$ and $\deg s = 2$. Put 
$f_i$ for the component of total degree $i$ in the power series expansion of $\log(1 + s +t)$. Then the map $\varphi_n:\R[s,t] \to \Val^{U(n)}$ of graded $\R$-algebras determined by
\begin{align}
\varphi_n(t) &\ge 0, \\
\label{eq:t^2 normalization}
		\varphi_n(t^2)& = \frac{2(2n-1)}{\pi}V^{2n}_2 , \\
\label{eq:s normalization}
		\varphi_n(s) & =\frac{n}{\pi} W^n_1
\end{align}
covers an isomorphism 
\begin{equation}
\C[s,t]/(f_{n+1},f_{n+2}) \simeq \Val^{U(n)}.
\end{equation}
\item The polynomials $f_i$ satisfy the relation
\begin{equation}\label{eq:relation relation}
nsf_n + (n+1)tf_{n+1} + (n+2) f_{n+2} = 0, \ n \ge 1
\end{equation}
and the diagram
\begin{equation}
\begin{CD}
\dots & & \dots \\
@VVV & & @V r VV \\
\C[s,t]/(f_{n+1},f_{n+2})& @> {\varphi_n}>> & \Val^{U(n)} \\
@VVV & & @V r VV \\
\C[s,t]/(f_{n},f_{n+1}) & @> {\varphi_{n-1}}>> &\Val^{U(n-1)} \\
\dots & & \dots \\
@VVV & & @V r VV \\
\C[s,t]/(f_{2},f_{3}) & @> {\varphi_1}>> &\Val^{U(1)} \\
\end{CD}
\end{equation}
commutes, where the vertical maps on the right are given by  restriction.
\end{itemize}
\end{theorem}

{\bf Remark.} Since
$$ \log(1 + s + t) =  (s + t) - \frac 1 2 (s+t)^2 + \frac 1 3 (s + t)^3 + \dots $$
it is easy to write down in closed form as many $f_i$ as desired. For example,
\begin{align*}
f_1 &= t ,\\
f_2 &= s - \frac 1 2 t^2, \\
f_3&= -st + \frac 1 3 t^3,\\
f_4 &= -\frac 1 2 s^2 + st^2 - \frac 1 4 t^4, 
\end{align*}
etc.

The rest of this paper is devoted to the proof of Theorem \ref{thm:main}. The reason for the coefficients in (\ref{eq:t^2 normalization}) and (\ref{eq:s normalization}) is the following.

\begin{lemma}\label{lm:values}
For $1\le k \le n$, put $r^n_k:\Val^{U(n)}\to\Val^{U(k)}$ to be the restriction map. Then
\begin{align}
		r^n_k\left( \frac{2(2n-1)}{\pi}V^{2n}_2\right) & =  \frac{2(2k-1)}{\pi}V^{2k}_2,\\
		r^n_k\left(\frac n \pi W^{n}_1\right) &=\frac k \pi W^k_1.
\end{align}
\end{lemma}
\begin{proof} It is clear that in each relation above the left and right sides are constant multiples of each other, so it will be enough to show that
\begin{align}
		r^n_1\left(\frac{2(2n-1)}{\pi} V^{2n}_2\right) & = \frac  2 \pi V^{2}_2,\\
		r^n_1\left(\frac n \pi W^{n}_1\right) &=\frac 1 \pi W^1_1,
\end{align}
 i.e. that 
 \begin{align}
 	\label{eq:t^2 values}
		 V^{2n}_2 (D_\C^1) & = \frac\pi{2n-1},\\
	\label{eq:s values}
		 W^{n}_1 (D_\C^1)&= \frac\pi n .
\end{align}

To prove (\ref{eq:t^2 values}) we note that the Hadwiger valuation $ V^{2n}_2$ of a smooth body $K\subset \C^n \simeq \R^{2n}$ may be expressed as the integral of the $(2n-3)$rd symmetric function of the principal curvatures of the boundary of $K$, multiplied by a certain constant $c$. Therefore
\begin{align}
\pi =V^{2n}_2(D_\C^n) &= c \binom{2n-1}{2n-3}\alpha_{2n-1}\\
&=  c \binom{2n-1}{2n-3} (2n) \frac{\pi^n}{n!},
\end{align}
so
\begin{equation} c= \frac{(n-2)!}{(4n-2)\pi^{n-1}}.
\end{equation}

On the other hand, the disk $D_\C^1$ is the Hausdorff limit of its tubular neighborhoods of radius $r$ as $r \downarrow 0$, which may be thought of as the union of $D_C^1 \times D_\C^{n-1}(r)$ together with a bundle of half-balls of real dimension $2n$ over the boundary circle. Computing the curvature integrals and passing to $r =0$, the boundary term tends to 0 and we obtain
\begin{align}
 V^{2n}_2(D_\C^1) = c \alpha_{2n-3} \pi = \frac{(n-2)!}{(4n-2)\pi^{n-1}} (2n-2) \frac{\pi^{n-1}}{(n-1)!}\pi = \frac\pi{2n-1},
\end{align}
as claimed.

To prove (\ref{eq:s values}) we apply Howard's transfer principle for Poincar\'e-Crofton formulas \cite{Howard:1993} to compare the integral geometry of $\C^n$ under the holomorphic isometry group $\overline{U(n)}$ with that of $\PP^n$ under its full isometry group $U(n+1)/U(1)$. To remain consistent with the measure on $\overline{U(n)}$ given in (\ref{eq:normalization}), we select the Haar measure on $U(n+1)/U(1)$ so that its total mass is equal to the volume $\frac {\pi^n}{n!}$ of $\PP^n$. The transfer principle then implies that 
\begin{align}\label{eq:transfer}
\frac 1 n = \frac{\int_{U(n+1)/U(1)} \#(\PP^1 \cap h \PP^{n-1}) \, dh}{\operatorname{area}(\PP^1)\operatorname{vol}(\PP^{n-1})}= 
\frac{\int_{\overline{U(n)}} \#(D_\C^1 \cap \bar g D_\C^{n-1}(R))\, d\bar g}
{\operatorname{area}(D_\C^1)\operatorname{vol}(D_\C^{n-1}(R))}
\end{align}

On the other hand (\ref{eq:W disks}) may be written
\begin{align}\label{eq:approx}
W^n_1(K)\sim \frac{\int_{\overline{U(n)}} \chi(K\cap \bar g D_\C^{n-1}(R))\, d\bar g }
{\operatorname{vol}(D_\C^{n-1}(R))}\end{align}
as $R \to \infty$. Applying this formula to $K = D_\C^1$,
(\ref{eq:transfer}) implies that $W^n_1(D_\C^1) = \frac \pi n$. 
\end{proof}

\begin{definition}\label{def:abuse} In view of this fact, for the sake of simplicity we will abuse notation by writing
$s$ for $ \frac{n}{\pi} W^n_1$,  and $t$ for the positive square root of
$  \frac{4n-2}{\pi}V^{2n}_2$. In computations with these elements the dimension in which we work should be clear from the context.
\end{definition}

\begin{corollary}\label{cor:volume constant} In $\Val^{U(n)}$, 
 $$t^{2n} = \frac{2\cdot(2n-1)!}{\pi^n(n-1)!}\vol^{2n}. $$
\end{corollary}
\begin{proof} We again use Howard's transfer principle, this time for the associated pairs $(\Rn,\overline{SO(n)})$ and $(S^n,SO(n+1))$. 

Put 
\begin{equation}
\Psi_i= \Psi^n_i:= \alpha_i^{-1}\alpha_{n-i-1}^{-1} \Phi_i^n.
\end{equation}
Thus the restriction of $\Psi_i$ to subsets of $\Rn$ of dimension $i$ is a multiple of the $i$-dimensional Hausdorff measure that transfers to $S^n$ to give
$$\Psi_i(S^i) = 1.$$
Now the kinematic formula for $\Rn$ may be expressed
$$c_n(\Psi_k) = \sum_{i+j =n+k} a_i \Psi_i \ot \Psi_j.$$
By the transfer principle the same formula applies to subsets of $S^n$--- this is true because the kinematic formula specializes to the Poincar\'e-Crofton formulas 
$$ \int \vol^k(M^i\cap \bar g N^j)\, d \bar g  = \int \Psi_k (M^i\cap \bar g N^j)\, d \bar g = c\Psi_i(M^i)\Psi_j(N^j) $$
for submanifolds $M^i,N^j$ with $i+j = n+k$. Taking the total measure of $SO(n+1)$ to be $\alpha_n$, and applying the formula to $S^i, S^j \subset S^n$, we find that all of the constants are equal to $\alpha_n$:
$$c_n(\Psi_k) =\alpha_n\sum_{i+j =n+k} \Psi_i \ot \Psi_j.$$
By Lemma \ref{lm:sesqui}, since $\Psi_0 = \frac 1 2 \chi = \frac 1 2$,
\begin{equation}
c_n(\Psi_k) = 2 c_n (\Psi_k\cdot \Psi_0) = 2 {\alpha_n}\sum_{i+j =n} (\Psi_k\cdot\Psi_i) \ot \Psi_j.
\end{equation}
Therefore $ \Psi_i \cdot \Psi_k = \frac 1 2 \Psi_{i+k}$, so $(4\Psi_1)^2 = 8\Psi_2$, and since 
\begin{align*}
8 \Psi_2(S^2)&= 8 ,\\
t^2 (S^2)& = \frac 2 \pi \area(S^2) = 8,
\end{align*}
it follows that 
\begin{equation}
t= 4\Psi_1
\end{equation}
and
\begin{equation}
t^n= 2^{n+1} \Psi_n = \frac{2^{n+1}}{\alpha_n }\vol^n. 
\end{equation}
In particular
\begin{equation}
t^{2n} = 2^{2n+1}\left(\frac{(2n-1)(2n-3)\dots 3 \cdot 1}{\pi^n 2^{n+1}}\right)\vol^{2n},
\end{equation}
as claimed.
\end{proof}

We note for future reference that
\begin{align}\label{eq:closed form for f}
f_k &= (-1)^{k+1}\sum_{i= 0}^{\left[ \frac {k} 2 \right]} \frac {(-1)^{i}} {k-i} \binom {k-i} i s^i t^{k-2i}\\
\notag &=(-1)^{k+1}\sum_{i= 0}^{\left[ \frac {k} 2 \right]} \frac {(-1)^{i}} {k-2i} \binom {k-i-1} i s^i t^{k-2i}.
\end{align}
In particular
\begin{equation}
f_{n+1}= (-1)^n \sum_{i=0}^{\left[ \frac {n+1} 2 \right]}
{(-1)^i}\frac{(n-i)!}{i!(n-2i+1)!}
s^i t^{n-2i+1}.
\end{equation}

\subsection{First deductions}
Our starting point is the following result of Alesker \cite{Alesker:2003}:
\begin{theorem}\label{thm:basis} The valuations
$$ U_{k,p}^n (K):= \int_{\bar G_\C(n,n-p)} t^{k-2p}(K\cap \bar P) \, d\bar P,$$
$0\le p \le \frac 1 2 \min\{k, 2n-k\}$, constitute a basis for $\Val^{U(n)}$. In particular, the Poincar\'e series of $\Val^{U(n)}$ is
\begin{equation}\label{eq:poincare series}
P_{\Val^{U(n)}}(x) = \frac{(1-x^{n+1})(1-x^{n+2})}{(1-x)(1-x^2)}.
\end{equation}
\end{theorem}
\begin{proof}
The first assertion is due  to Alesker in \cite{Alesker:2003}. The assertion about the Poincar\'e series then follows from a simple comparison between the coefficients of the given polynomial and Alesker's computation of the dimensions of the $\Val^{U(n)}_k$.
\end{proof}

As in the discussion preceding the statement of Theorem \ref{thm:valSOn}, these valuations are monomials in $s$ and $t$:

\begin{proposition}\label{prop:monomial basis}
$$U_{k,p}^n = s^p t^{k-2p}.\hskip 1in \square$$
\end{proposition}

This establishes the first assertion of Theorem \ref{thm:main}.

{\bf Remark.} Thus the basis described in Theorem \ref{thm:basis} may be understood as follows. Up through degree $i=n$, there are no relations between $s$ and $t$. From degree $i=n+1$ through the highest degree $i=2n$, the basis elements are in one-to-one correspondence with those in degree $2n-i$: in fact they are simply the products of the latter with $t^{2i-2n}$. 
 
From general considerations we also find 

\begin{lemma} There are polynomials $p_{n+1},p_{n+2}$, of degrees $n+1, n+2$ respectively, such that
\begin{equation}
\ker \varphi_n = (p_{n+1}, p_{n+2}).
\end{equation}
\end{lemma}
\begin{proof} By Alesker's basis Theorem \ref{thm:basis} there are relations $p_{n+1}, p_{n+2}\in \C[s,t]$ in the given degrees such that $p_{n+2} \ne t\cdot p_{n+1}$. We first show that these polynomials are relatively prime.  Otherwise, let $ w \in \C[s,t]$ be an element of degree $0<k < n+1$ dividing both, with $p_j = w d_j, j = n+1,n+2$, where the $d_j$ are relatively prime. It is clear that all of these elements are homogeneous. Therefore $W:=\C[s,t]/(d_{n+1}, d_{n+2})$ is a graded algebra with Poincar\'e series
\begin{equation}
P_W(x) = \frac{(1-x^{n+1-k})(1-x^{n+2-k})}{(1-x)(1-x^2)},
\end{equation}
which has degree $2n-2k$. Thus the image of the linear map $h:W \to \Val^{U(n)}$ covered by multiplication by $w$ in $\C[s,t]$ meets the socle $\Val^{U(n)}_{2n}$ only at $0$. However, this contradicts the fact that multiplication induces a perfect pairing on $\Val^{U(n)}$: since the image of $w$ in $\Val^{U(n)}$ is nonzero, there exists $g \in \C[s,t]$, $\deg g = 2n-k$, such that $w \cdot g \ne 0$ in $\Val^{U(n)}$.  Thus the image of $g$ in $W$ under $h$ is not zero. This is a contradiction.

It follows that the Poincar\'e series of $\C[s,t]/(p_{n+1},p_{n+2})$ is given by (\ref{eq:poincare series}). Since $\Val^{U(n)}$ is isomorphic to a quotient of this algebra we obtain the desired conclusion.
\end{proof}

From this and Lemma \ref{lm:values} we find that the last conclusion of Theorem \ref{thm:main} is valid with the $f_i$ replaced by the polynomials $p_i$, to be determined. Thus only the second assertion remains. The combinatorial part (\ref{eq:relation relation}) follows at once by writing
\begin{align}
\log(1 + tx + sx^2)& = \sum_{i= 1} ^\infty f_i(s,t) x^i, \\
f_1 + (tf_1 +2f_2) x +\sum_{i= 1} ^\infty [is f_i + (i+1)tf_{i+1} + (i+2)f_{i+2}] x^{i+1}&=(1 + tx + sx^2)\sum_{i= 1} ^\infty i f_i x^{i-1} \\
& =(1 + tx + sx^2)\frac d {dx}  \log(1 + tx + sx^2) \\
& = t + 2sx.
\end{align}

\subsection{Identifying $\mathcal A^n$}

For the rest of the paper we will abbreviate $\mathcal A^n:= \mathcal A^{U(n)}$. Using Corollary \ref{cor:annihilator} we can give an explicit basis for this subspace. 

\begin{lemma} \label{lm:A basis}
For $2\le j \le 2n-2$, the elements
\begin{equation}\label{eq: A basis}
(n-i)s^it^{j-2i} - (4n-4i -2)s^{i+1}t^{j-2i-2 }, \ 0\le i \le \frac{\min\{j, 2n-j\}}  2- 1
\end{equation}
constitute a basis for $\mathcal A^n_j$.
\end{lemma}
\begin{proof}
By Theorem \ref{thm:basis}, together with $t^j$ the elements (\ref{eq: A basis}) constitute a basis for $\Val^{U(n)}_j$.  Hence it is enough to show that they belong to $\mathcal A^n$. 

Put $r_k:= k t^2- (4k-2)s$. By Lemma \ref{lm:values}, $r_k \in \mathcal A^k$. By the definition of the generating valuation $s$, it follows that $s^{n-k}r_k \in \mathcal A^n$ for $k\le n$. Now the relation (\ref{eq:inclusion}) implies that $t^is^{n-k}r_k \in \mathcal A^n$ for all $i \ge 0$. These include the elements (\ref{eq: A basis}).
\end{proof}

Denote the ordered monomial basis $(t^j, st^{j-2},\dots)$ for $\Val^{U(n)}_j$ given in Theorem \ref{thm:basis} and Prop. \ref{prop:monomial basis} by $\mathfrak b_j$, and the ordered basis consisting of $t^j$ together with the degree $j$ elements of (\ref{eq: A basis}) by $\mathfrak c_j$. Thus if  we define the $(k+1) \times (k+1)$ matrix
\begin{equation}\label{eq:Ank}
A_k^n:= \left [\begin{matrix} 	1  & 0 & 0 & 0 & \dots & 0 & 0\\
					n & -2(2n-1) & 0 & 0 &\dots & 0 & 0\\
					0 & (n-1) &-2( 2n -3) & 0 & \dots & 0 & 0\\
					\dots \\
					0 & 0 & 0 & 0 & \dots & n -k + 1 & -2(2n-2k + 1)
		\end{matrix}\right]
\end{equation}
then for $2k+1 \le n$ we have
\begin{equation}
\mathfrak c_j  = A^n_k \mathfrak b_j, 
\end{equation}
$  j = 2k, 2k +1, 2n-2k-1,2n-2k$.

\subsection{Framework for induction}
Recalling the remark following Prop. \ref{prop:monomial basis}, if $2k + 1 \le n$ then the $\R$ vector spaces
$$\Val^{U(n)}_{2k}, \Val^{U(n)}_{2k+1},\Val^{U(n)}_{2n-2k-1},\Val^{U(n)}_{2n-2k} $$
all have dimension $k+1$. In fact the maps
\begin{equation}
 a \mapsto t^{j}\cdot a,
\end{equation}
$j = 1, 2n-4k-1, 2n-4k$, 
are isomorphisms $\Val^{U(n)}_{2k} \to \Val^{U(n)}_{2k+1},\Val^{U(n)}_{2n-2k-1},\Val^{U(n)}_{2n-2k} $ respectively.
Since the Poincar\'e duality pairing is given by multiplication it is trivial to see:
\begin{proposition}  The pairings on $\Val^{U(n)}_{2k}$  given by
\begin{equation}
\langle a,b\rangle:= \PD(a, t^{2n-4k}b)
\end{equation}
and
\begin{equation}
\langle\langle a,b\rangle\rangle:= \PD(ta, t^{2n-4k-1}b)
\end{equation}
are identical to one another, and are symmetric and nondegenerate. $\square$
\end{proposition}

We denote by $P^n_k$ the $(k+1) \times (k+1)$ matrix giving this pairing with respect to the ordered monomial basis $t^{2k}, st^{2k-2}, \dots, s^k$, and put $Q^n_k: = \left(P^n_k\right)^{-1}$. Thus $Q^n_k$ is the matrix giving the associated pairing on the dual spaces. By Corollary \ref{thm:KG and PD}, $Q^n_k$ is also the matrix of coefficients for the kinematic formula $k_{n}(1)$.

{\bf Convention.} From this point on we normalize the $k_n$ so that equation (\ref{eq:annihilator}) in Prop. \ref{prop:annihilator} is literally true.

With this convention, Corollary \ref{cor:volume constant} gives
\begin{align}
\notag \int_{\overline{U(n)}} \chi(K \cap \bar g L) \, d\bar g &= \chi(K) \vol(L) + \dots\\
\notag &=\frac{\pi^n(n-1)!}{2\cdot(2n-1)!}( 1 \ot t^{2n} +\dots)(K,L) \\
\label{eq:refined kinematic}&= \frac{\pi^n(n-1)!}{2\cdot(2n-1)!}k_{n}(1)(K,L) \\
\end{align}
where
\begin{align}
k_n(1) =\sum_{i= 0}^{2n}
(t^i, s t^{i-2},\dots, s^k t^{i-2k}) \ot Q^n_{\min\left\{\left[\frac i 2\right],\left[\frac{2n-i}2\right]\right\}}\left(\begin{matrix}
t^{2n-i}\\
s t^{2n-i-2}\\
\dots \\
s^kt^{2n-i-2k}
			\end{matrix}\right).
\end{align}
Here $Q^n_0 = 1.$ Note that in this last sum there are four terms involving each of the matrices $Q^n_k$, $k\le \frac n 2 -1$; on the other hand there are three such terms when $n$ is odd and $k= \frac{n-1} 2$, and one such term when $n$ is even and $k= \frac n 2$.

Now by Prop.  \ref{prop:annihilator}, for $k\ge 1$ there exists a nonsingular symmetric $k\times k$ matrix $\widetilde{Q^n_k}$ such that
\begin{equation}\label{eq:Qtilde}
Q^n_k = (A^n_k)^t \left[
					\begin{matrix}
								1 & \vec 0 \\
								\vec 0 & \widetilde{Q^n_k}
							\end{matrix}
							\right] A^n_k
\end{equation}

We will determine the relations in $\Val^{U(n)}$ between the generators $s$ and $t$ using an induction on $n $ and $k$. 

\begin{definition}
For each $\bar P \in \bar G_\C(n+1,n)$, choose a holomorphic isometry $\gamma_{\bar P}:\bar P \to \C^n$, and define the map $\iota: \Val^{U(n)} \to \Val^{U(n+1)}$  by
$$\iota(\mu)(K) := \int_{\overline G_\C(n+1,n)} \mu(\gamma_{\bar P}(K\cap \bar P)) \, d\bar P.$$
Clearly this is independent of the choices  of the $\gamma_{\bar P}$.
\end{definition}

\begin{proposition}\label{prop:step-up}
Under the maps
\begin{equation}
\begin{CD}
\Val^{U(n+1)}\ot \Val^{U(n+1)} @> id \ot r>>\Val^{U(n+1)}\ot \Val^{U(n)} @< \iota\ot id<<  \Val^{U(n)}\ot \Val^{U(n)},
\end{CD}
\end{equation}
the elements $k_{n+1}(1)$ and $k_n(1)$ satisfy
\begin{equation}\label{eq:step-up}
\frac \pi{2(2n+1)}(id \ot r) (k_{n+1}(1) )= (\iota \ot id)(k_n(1)).
\end{equation}
In terms of the notational abuse of Definition \ref{def:abuse}, this relation may be written 
\begin{equation}\label{eq:stepup 2}
(n+1) k_{n+1}(1) = 2(2n+1)(s\ot 1)\cdot k_n(1)
\end{equation}
as elements of $\Val^{U(n+1)}\ot \Val^{U(n)}$.
\end{proposition}

\begin{proof} Clearly the images of $s$ and $t$ under the restriction map satisfy
\begin{equation}
r(s_{n+1} ) = s_n, \ r(t_{n+1}) = t_n.
\end{equation}
Furthermore, in view of our normalizing conventions for the various Haar measures,
\begin{equation}\label{eq:times s}
\iota(s_n^it_n^j) = \frac \pi {n+1} s_{n+1}^{i+1}t_{n+1}^j.
\end{equation}
Returning to our usual abuse of notation we write more succinctly
\begin{align}
r(s ) = s, &\ r(t) = t, \\
\label{eq:times s2}\iota(s^it^j)& = \frac \pi {n+1} s^{i+1}t^j.
\end{align}


To prove the relation (\ref{eq:step-up}),  for $\bar P \in \bar G_\C(n+1,n)$ we put $\overline{U(n)}_{\bar P}$ for the left coset  of $\overline{U(n)}$ in $\overline{U(n+1)}$ consisting of the elements that map $\C^n$ to $\bar P$. Now if $K \in \mathcal K(\C^{n+1})$ and $L \in \mathcal K(\C^n) \subset \mathcal K(\C^{n+1})$ then (denoting by $i: \C^n \to \C^{n+1}$ the inclusion)
\begin{align*}
\frac{\pi^{n+1} n!}{2\cdot (2n+1)!}(id \ot r)(k_{n+1}(1))(K,L) &=\frac{\pi^{n+1} n!}{2\cdot (2n+1)!}k_{n+1}(1) (K,i(L)) \\
&= \int_{\overline{U(n+1)}} \chi (K \cap \bar g L) \, d\bar g \\
& = \int_{\overline{G}_\C(n+1,n)}\int_{\overline{U(n)}_{\bar P}}\chi (K \cap \bar g L) \, d\bar g \, d \bar P \\
& =\int_{\overline{G}_\C(n+1,n)}\int_{\overline{U(n)}_{\bar P}}\chi ((K \cap\bar P) \cap \bar g L) \, d\bar g \, d \bar P \\
&= \int_{\overline{G}_\C(n+1,n)}\int_{\overline{U(n)}}\chi (\gamma_P(K \cap\bar P) \cap \bar h L) \, d\bar h \, d \bar P \\
&= \frac{\pi^{n} (n-1)!}{2\cdot (2n-1)!}\int_{\overline{G}_\C(n+1,n)} k_n(1) (\gamma_{\bar P}(K\cap \bar P), L)
 \, d \bar P \\
&= \frac{\pi^{n} (n-1)!}{2\cdot (2n-1)!}(\iota \ot id)( k_n(1))(K,L),
\end{align*}
which simplifies to (\ref{eq:step-up}). The relation (\ref{eq:stepup 2}) follows from this and (\ref{eq:times s2}).
\end{proof}

\subsection{Proof of theorem}

\begin{lemma}\label{lm:PQ}
If  $2k\le n-1$ then 
\begin{equation}\label{eq:PQ matrix}
\frac {n}{2(2n-1)}Q^{n}_{k}P^{n-1}_{k} = 
\left(\begin{matrix}
0 & 0 & 0 & 0 &\dots & -a_0^{n,k} \\
1& 0 & 0 & 0 & \dots & -a_1^{n,k} \\
0 & 1 & 0 & . & \dots & -a_2^{n,k} \\
 0& 0 & 1 & . & \dots & -a_3^{n,k} \\
 & . & . & . & \dots & .  \\
. & . & . & . & \dots & . \\
0 & 0 & 0 & \dots & 1 &-a_{k}^{n,k}\\
\end{matrix}\right),
\end{equation}
where, putting
\begin{align}
\notag a_{k+1}^{n,k} &:= 1,\\
\varphi^{n,k} &:= \sum_{i=0}^{k+1} a_i^{n,k} s^i t^{2n-2k-2i-1},
\end{align}
we have 
\begin{equation}\label{eq:main relation}
\varphi^{n,k} = 0 \text{ in } \Val^{U(n)},  \ 0\le k \le \frac n 2 - 1;
\end{equation}
and if $n$ is odd then 
\begin{equation}
t \varphi^{n, \frac{n-1}{2}} =0 \text{ in }\Val^{U(n)}.
\end{equation}
\end{lemma}

Observe that in the latter case $t\varphi^{n,\frac {n-1} 2}$ is a polynomial even though $\varphi^{n,\frac {n-1} 2}$ is only rational.

\begin{proof} By (\ref{eq:refined kinematic}), the terms of bidegree $(2n-2k, 2k)$ in $k_n(1)$ may be written
\begin{equation}\label{eq:val^n relation}
\left(t^{2n-2k},\dots, s^kt^{2n-4k}\right) \ot Q^n_k \left(\begin{matrix}
t^{2k}\\
\dots \\
s^k
\end{matrix}\right) \in \Val^{U(n)}\ot\Val^{U(n)},
\end{equation}
and similarly the terms of bidegree $(2n-2k-2, 2k)$ in $k_{n-1}(1)$ are given by
\begin{equation}\label{eq:val^n-1 relation}
\left(t^{2n-2k-2},\dots, s^kt^{2n-4k-2}\right) \ot Q^{n-1}_k \left(\begin{matrix}
t^{2k}\\
\dots \\
s^k
\end{matrix}\right) \in \Val^{U(n-1)}\ot\Val^{U(n-1)}.
\end{equation}
Therefore the relation (\ref{eq:step-up}) gives
\begin{equation}
\frac{n}{2(2n-1)}\left(t^{2n-2k},\dots, s^kt^{2n-4k}\right) \ot Q^{n}_k \left(\begin{matrix}
t^{2k}\\
\dots \\
s^k
\end{matrix}\right) = 
\left(st^{2n-2k-2},\dots, s^{k+1}t^{2n-4k-2}\right)\ot Q^{n-1}_k \left(\begin{matrix}
t^{2k}\\
\dots \\
s^k
\end{matrix}\right)\end{equation}
in $\Val^{U(n)}\ot\Val^{U(n-1)}$. Since $t^{2k},\dots, s^k$ constitute a basis for $\Val^{U(n-1)}_{2k}$ it follows that
\begin{equation}
\frac{n}{2(2n-1)}\left(t^{2n-2k},\dots, s^kt^{2n-4k}\right) Q^{n}_k  = 
\left(st^{2n-2k-2},\dots, s^{k+1}t^{2n-4k-2}\right) Q^{n-1}_k, \end{equation}
which shows that the left hand side of (\ref{eq:PQ matrix}) has the given form for some constants $a^{n,k}_i$, where the resulting polynomial $\varphi^{n,k}$ satisfies $t \varphi^{n,k}=0$ in $\Val^{U(n)}$.

If $k \le \frac n 2 - 1$ then the relations corresponding to (\ref{eq:val^n relation}) and (\ref{eq:val^n-1 relation}) hold for the terms of bidegree $(2n-2k-1,2k+1)$ and $(2n-2k-3,2k+1)$. Arguing as above we then arrive at the relation (\ref{eq:main relation}).
\end{proof}

To prove theorem we will show that if $n$ is odd (resp. even) then $t\varphi^{n,\frac{n-1} 2}$(resp. $\varphi^{n, \frac n 2-1}$ ) is a nonzero constant multiple of $f_{n+1}$. Since the natural map $\Val^{U(n+1)} \to \Val^{U(n)}$  is well-defined and $f_{n+2} = f_{(n+1)+1}=0$ in $\Val^{U(n+1)}$, it follows that $f_{n+2} = 0 $ in $\Val^{U(n)}$ as well.

In fact we prove more generally:

\begin{proposition}\label{prop:general term} For $ 2k \le n-1$ and $0\le i\le  k$,
\begin{equation}\label{eq:general term}
a_i^{n,k} = (-2)^{i-k-1} \binom {k+1}{i} \frac{(n-i)(n-i-1)\dots(n-k)}{(2n-2k-2i-1)(2n-2k-2i-3)\dots (2n-4k-1)} .
\end{equation}
\end{proposition} 

\begin{corollary}
If $n$ is even then
\begin{equation} \varphi^{n,\frac n 2-1} = (-1)^{\frac n 2 } f_{n+1}.
\end{equation}
If n is odd then
\begin{align} 
t\varphi^{n,\frac{n-1} 2}&= (-1)^{\frac {n-1} 2 } \left(\frac {n+1} 2\right) f_{n+1}, \\
\end{align}
\end{corollary}

The rest of this section is devoted to the proof of Prop. \ref{prop:general term} using induction on $k$, the starting point being the relation
\begin{equation}\label{eq:start induction}
a^{n,0}_0 = \frac{-n}{2(2n-1)}, \ n \ge 2,
\end{equation}
which is clearly valid from the defining relation (\ref{eq:PQ matrix}) since $Q^n_0 = 1$ for all $n$.

\begin{lemma} \label{lm:induction} 
Let
$$ R^n_k := \frac n{2(2n-1)}\left[ \begin{matrix} 
\frac{2(2n-1)}{n} & 1 & \frac{2\binom{2n-5}{n-2}}{\binom{2n-2}{n-1}} & \frac{ 2 \binom{2n-7}{n-3}}{\binom{2n-2}{n-1}}& \dots &  \frac{2 \binom{2n-2k+1}{n-k+1}}{\binom{2n-2}{n-1}}& \frac{2 \binom{2n-2k-1}{n-k}}{\binom{2n-2}{n-1}}\\
1       & 0 &            0           &                  0              & \dots & 0 & -a^{n-1,k-1}_0\\
0      & 1 &            0           &                  0              & \dots & 0 & -a^{n-1,k-1}_1\\
\dots      & \dots&            \dots           &                  \dots              & \dots & \dots & \dots\\
0       & 0 &            0           &                  0              & \dots & 1 & -a^{n-1,k-1}_{k-1}\\
\end{matrix}\right] .$$
Then
\begin{equation}
 R^n_kQ^n_k \\
=  \left[\begin{matrix}
1 & 0 \\
0 & Q^{n-1}_{k-1}
\end{matrix}\right].
\end{equation}
\end{lemma}

{\bf Remark.} Note that the first row of $R^n_k$ may also be written as
\begin{equation}\label{eq:first row}
\binom{2n-1} n^{-1} \left({\binom{2n-1}{n}}, {\binom{2n-3}{n-1}},\dots, {\binom{2n-2k-1}{n-k}} \right)
\end{equation}

\begin{proof} Equating the terms of bidegree $(2n-2k,2k)$ in Prop. \ref{prop:step-up}, we obtain for $0 \le 2k\le n$:
\begin{align*}
\frac n {2(2n-1)}(t^{2n-2k}, st^{2n-2k-2},\dots,s^k t^{2n-4k}) \ot Q^n_k 
\left(\begin{matrix}
t^{2k} \\
st^{2k-2}\\
.\\
.\\
.\\
s^k
\end{matrix}\right) &= \\
=(t^{2n-2k}, st^{2n-2k-2},&\dots,s^{k-1} t^{2n-4k+2})\ot Q^{n-1}_{k-1} 
\left(\begin{matrix}
st^{2k-2}\\
.\\
.\\
.\\
s^k
\end{matrix}\right)
\end{align*}
in $\Val^{U(n-1)}\ot\Val^{U(n)} $.
Applying Lemma \ref{lm:PQ}, it follows that
\begin{equation}
(t^{2n-2k}, st^{2n-2k-2},\dots,s^k t^{2n-4k}) = 
(t^{2n-2k}, st^{2n-2k-2},\dots,s^{k-1} t^{2n-4k+2}) 
\left[I_k | -\vec a^{n-1,k-1} \right]
\end{equation}
in $\Val^{U(n-1)}$. Since $t^{2n-2k}, st^{2n-2k-2},\dots,s^{k-1} t^{2n-4k+2}$ constitute a basis for $\Val^{U(n-1)}_{2n-2k}$, and $t^{2k},\dots, s^k$ are a basis for $\Val^{U(n)}_{2k}$, recalling (\ref{eq:Qtilde}) we find that
\begin{align}
\notag \frac n {2(2n-1)} \left[I_k | -\vec a^{n-1,k-1}\right](A^n_k)^t \left[
					\begin{matrix}
								1 & \vec 0 \\
								\vec 0 & \widetilde{Q^n_k}
							\end{matrix}
							\right] A^n_k
&=\frac n {2(2n-1)} \left[I_k | -\vec a^{n-1,k-1}\right] Q^n_k \\
\label{eq:second} &= \left[\vec 0 \ | \ Q^{n-1}_{k-1}\right].
\end{align}
In view of the definition (\ref{eq:Ank}) of $A^n_k$, the first row of $ \left[\begin{matrix}
								1 & \vec 0 \\
								\vec 0 & \widetilde{Q^n_k}
							\end{matrix}
							\right] A^n_k$
is $e_1:=(1,0,\dots,0)$. Now if we denote by $M$ the $(k+1)\times (k+1)$ matrix whose first row is $e_1$ and whose bottom $k$ rows are identical to  $\frac n {2(2n-1)} \left[I_k | -\vec a^{n-1,k-1}\right](A^n_k)^t$ then (\ref{eq:second}) becomes
\begin{equation}\label{eq:temp}
M \left[\begin{matrix}
								1 & \vec 0 \\
								\vec 0 & \widetilde{Q^n_k}
							\end{matrix}
							\right] A^n_k = 
							\left[\begin{matrix}
								1 & \vec 0 \\
								\vec 0 & {Q^{n-1}_{k-1}}
							\end{matrix}
							\right].
\end{equation}
Using the identity
$$
(n-i)\binom {2n-2i-1}{n-i}- 2(2n-2i-1) \binom{2n-2i-3}{n-i-1}=0,$$
 a straightforward calculation  shows that
\begin{equation*}
M = R^n_k (A^n_k)^t
\end{equation*}
so the desired relation follows from (\ref{eq:Qtilde}) and (\ref{eq:temp}).
\end{proof}

To complete the proof of theorem we use Lemma \ref{lm:induction} to write
\begin{align*}
 Q^n_k P^{n-1}_k &=
\left(R^n_k\right)^{-1} 
\left[\begin{matrix}
								1 & \vec 0 \\
								\vec 0 & {Q^{n-1}_{k-1}}
							\end{matrix}
							\right] 
\left[\begin{matrix}
								1 & \vec 0 \\
								\vec 0 &{P^{n-2}_{k-1}}
							\end{matrix}
							\right] 
							R^{n-1}_k \\
&= \left(R^n_k\right)^{-1} 
\left[\begin{matrix}
								1 & \vec 0 \\
								\vec 0 & {Q^{n-1}_{k-1}}{P^{n-2}_{k-1}}
							\end{matrix}
							\right]  
							R^{n-1}_k,
\end{align*}
yielding
\begin{equation}
\frac{n-1}{2(2n-3)} R^n_k \left(\frac{n}{2(2n-1)}Q^n_kP^{n-1}_k\right)
= \frac{n}{2(2n- 1)}\left[
\begin{matrix}
\frac{n-1}{2(2n-3)} & \vec 0  \\
\vec 0  & \frac{n-1}{2(2n-3)}Q^{n-1}_{k-1}P^{n-2}_{k-1}  
\end{matrix}
\right] R^{n-1}_k.
\end{equation}
Recalling Lemma \ref{lm:PQ} and the definition of $R^n_k$, if we equate the last columns on the right and the left we obtain the following relations among the constants $a^{n,k}_i, a^{n-1,k-1}_i,a^{n-2,k-1}_i$:
\begin{align*}
-\sum_{i=0}^k \binom{2n-2i -1}{n-i} a^{n,k}_i& =
\binom{2n-2k-3}{n-k-1}, \\
-a^{n,k}_i + a^{n-1,k-1}_{i}a^{n,k}_k &=
-a_{i-1}^{n-2,k-1}+ a^{n-1,k-1}_{i}a^{n-2,k-1}_{k-1}, \ i =0,\dots,k-1,
\end{align*}
where we set $a^{n-2,k-1}_{-1} := 0$. Equivalently,
\begin{align}
\label{eq:a1}\sum_{i=0}^{k+1} \binom{2n-2i -1}{n-i} a^{n,k}_i & = 0 \\
\label{eq:a2} a^{n-2,k-1}_{i-1} + a^{n-1,k-1}_{i}\left(a^{n,k}_k- a^{n-2,k-1}_{k-1}\right)& =a_i^{n,k}  , \ i =0,\dots,k-1.
\end{align}
Since all of the matrices above are invertible these relations determine the $a_i^{n,k}$ uniquely in terms of the $a^{n-1,k-1}_i,a^{n-2,k-1}_i$, i.e. the system above is nonsingular in $a_0^{n,k},\dots,a_k^{n,k}$. Therefore, to complete the proof of Prop. \ref{prop:general term} by induction on $k$ it is enough to show that (\ref{eq:a1}) and (\ref{eq:a2}) are valid for the stated  values (\ref{eq:general term}). 

Starting with the observation that, with these values,
$$
a^{n,k}_k - a^{n-2,k-1}_{k-1} = \frac{-n}{2(2n-4k-1)},
$$
the verification of (\ref{eq:a2}) is a straightforward calculation. Meanwhile, fixing $n$ and $k$, the relation (\ref{eq:a1}) reduces to 
\begin{equation}\label{eq:reduction}
\sum_{i=0}^{k+1}  (-1)^i \binom{k+1}{i}\frac{(2n-2i-1)\dots(2n-2k-1)}{(2n-2k-2i-1)\dots(2n-4k-1)} = 0,
\end{equation}
where the $(k+1)$st term is understood to be $(-1)^{k+1}$--- in fact, after substituting the values (\ref{eq:general term}) into (\ref{eq:a1}) we find that the ratio of the $i$th terms of respective left-hand sides of (\ref{eq:a1}) and (\ref{eq:reduction}) is
\begin{equation}
\frac{(2n-2k-3)(2n-2k-5)\dots \cdot 3 \cdot 1 (-1)^{k+1} 2^{n-k-2}}{(n-k-1)!},
\end{equation}
independent of $i$.
Substituting $z:= \frac{2n-1}{2}$ and multiplying the sum by the function
$(z-k-1)(z-k-2)\dots(z-2k)$, the relation (\ref{eq:reduction}) becomes
\begin{equation}\label{eq:last relation}
\sum_{i=0}^{k+1}  (-1)^i \binom{k+1}{i}{(z-i)(z-i-1)\dots(z-i-k+1)} = 0.
\end{equation}
If $\Delta$ is the difference operator $\Delta(f(z)):= f(z) -f(z-1)$ then the left-hand side is $$\Delta^{k+1}\left(z(z-1)\dots(z-k+1)\right),$$ which vanishes identically since the subject polynomial has degree $k$. $\square$

\section{Open questions}

{\bf 1.} Obviously $\Val^{U(n)}$ is much more complicated than $\Val^{SO(n)}$, and many questions that are trivial in the latter case are not in the former. As we have seen, the Poincar\'e duality pairing for $\Val^{SO(n)}$ is essentially as simple as possible, so the deduction of the kinematic formulas via Thm. \ref{thm:KG and PD} is very easy. For $\Val^{U(n)}$ the question of determining the pairing matrices $P^n_k$ and their inverses $Q^n_k$, which determine the kinematic formulas, is open. Using, the MAGMA computer algebra package, Graham Matthews has calculated the $Q^n_k$ for $k\le 11$. The entries are rational functions of $n$, with numerators and denominators having irreducible factors of even degree apparently growing without bound as $k$ increases. We have not been able to discern the patterns in the coefficients.

A seemingly simpler question is: are the $Q^n_k$ positive definite? They are for small values of $n$. 

{\bf 2.} Hadwiger's basis theorem for $\Val^{SO(n)}$ implies that various methods for constructing $SO(n)$-invariant valuations lead to the same results. For example, if $K$ is a compact convex body then $t^i(K)$ is equal to the average value of $t^i(\pi_P(K))$ as $P$ ranges over the Grassmannian $G(n,j)$, $j \ge i$; or, if the body has smooth boundary, as the integral over the boundary of $K$ of the $(n-i-1)$st elementary symmetric functions of the principal curvatures. In the unitary case, we have seen that the monomials $s^it^j$ correspond to Alesker's $U_{k,p}$ basis, given by integrating the Hadwiger valuations of the intersections of $K$ with affine complex planes. On the other hand Alesker also defined an alternative basis, denoted $C_{k,p}$, given by averaging the Hadwiger valuations of the {\it projections} of $K$ to the elements of the various complex Grassmannians. Furthermore, H. Park \cite{Park:2002} has classified the $\overline{U(n)}$-invariant differential forms on the sphere bundle of $\C^n$, whose pairings with $N(K)$ give valuations in $\Val^{U(n)}$ as before. Determining the linear relations among these bases is a fundamental open problem.

{\bf 3.} Say that a valuation $\varphi$ is {\it positive} if $\varphi(K)\ge 0$ for all convex bodies $K$, and {\it monotone }if $\varphi(K)\ge \varphi(L)$ whenever $K\supset L$. It is easy to see that the cones of positive and monotone valuations in $\ValSO$ coincide, and consist of all nonnegative linear combinations of the Hadwiger valuations. What are the positive and monotone cones in $\ValU$?

\bibliographystyle{plain}
\bibliography{unitary-4-11}

\end{document}